%% file: ms.tex
\documentclass[review=false, nonacm=true, authorversion=true]{acmart}
\setcopyright{none}
\usepackage{graphicx}
\usepackage{subcaption}
\usepackage{amsmath}
\usepackage[utf8]{inputenc}
\usepackage{backnaur}
\usepackage{physics}
\usepackage{moresize}
\usepackage{algorithm}
\usepackage{algpseudocode}
\usepackage{tikz}
\usetikzlibrary{shapes}
\renewcommand\footnotetextcopyrightpermission[1]{} % removes footnote with conference information in first column
\title{Optimizing Geometric Multigrid Methods with Evolutionary Computation}
\author{Jonas Schmitt}
\orcid{0000-0002-8891-0046}
\email{jonas.schmitt@fau.de}
\author{Sebastian Kuckuk}
\email{sebastian.kuckuk@fau.de}
\author{Harald Köstler}
\email{harald.koestler@fau.de}
\affiliation{%
	\institution{Friedrich-Alexander University Erlangen-Nürnberg (FAU)}
	\department{Chair for System Simulation, Department of Computer Science}
	\streetaddress{Cauerstraße 11}
	\city{Erlangen}
	\postcode{91058}
	\country{Germany}
}

\begin{document}

\begin{abstract}
For many linear and nonlinear systems that arise from the discretization of partial differential equations the construction of an efficient multigrid solver is a challenging task.
Here we present a novel approach for the optimization of geometric multigrid methods that is based on evolutionary computation, a generic program optimization technique inspired by the principle of natural evolution.
A multigrid solver is represented as a tree of mathematical expressions which we generate based on a tailored grammar.
The quality of each solver is evaluated in terms of convergence and compute performance using automated local Fourier analysis (LFA) and roofline performance modeling, respectively.
Based on these objectives a multi-objective optimization is performed using strongly typed genetic programming with a non-dominated sorting based selection.
To evaluate the model-based prediction and to target concrete applications, scalable implementations of an evolved solver can be automatically generated with the ExaStencils framework.	
We demonstrate our approach by constructing multigrid solvers for the steady-state heat equation with constant and variable coefficients that consistently perform better than common V- and W-cycles.
\end{abstract}

\maketitle
\thispagestyle{empty}

\section{Introduction}\label{sec:intro}
\input{src/introduction}

\section{A Formal Grammar for Generating Multigrid Solvers}\label{sec:grammar}
\input{src/grammar}

\section{Multi-Objective Optimization with Evolutionary Computation}\label{sec:optimization}
\input{src/optimization}

\section{Evolving Solvers for the Steady-State Heat Equation}\label{sec:experiments}
\input{src/experiments}

\section{Conclusion}\label{sec:conclusion}
\input{src/conclusion}
\bibliographystyle{acm}
\bibliography{references}

\end{document}

%% file: src/introduction.tex
Solving the linear or nonlinear systems that arise from the discretization of partial differential equations (PDEs) efficiently is an outstanding challenge.
The huge number of unknowns in many of these systems necessitates the design of fast and scalable solvers.
Unfortunately, the optimal solution method highly depends on the system itself and it is therefore infeasible to formulate a single algorithm that works efficiently in all cases.
Multigrid methods are a class of asymptotically optimal solution algorithms for (non-)linear systems.  
Although in the last decades great effort has been put into the design of efficient multigrid methods for many important cases, such as the Navier-Stokes or Maxwell's equation, this task is still not fully solved.
Within this paper we propose a novel approach for the automatic optimization of multigrid solvers through the use of evolutionary computation.
Our approach builds on the work by Mahmoodabadi and Köstler \cite{Mahmoodabadi2017}, in which it was first demonstrated how genetic programming (GP) \cite{koza1994genetic} can be used to optimize the iteration matrices of stationary iterative methods.
In contrast to this work, where the iteration matrix was assembled explicitly, we represent iterative methods as symbolic expressions, which are independent of the size of the linear system.
For this purpose we propose a new formal grammar for the automatic generation of multigrid expressions. 
We furthermore show how we can automatically obtain convergence and performance estimates for geometric multigrid solvers on rectangular grids of arbitrary size and how, based on these metrics, a multi-objective optimization can be performed using genetic programming and evolution strategies (ES) \cite{beyer2002evolution}.
In addition the evolved multigrid solvers are emitted in the form of a domain specific language (DSL) specification and the ExaStencils code generation framework is employed to automatically generate a scalable implementation on the target platform.
Our approach is therefore similar to the work by Thekale et al \cite{thekale2010optimizing}, where the optimal number of full multigrid cycles was optimized based on a cost model, but aims to achieve more generality by considering the construction of multigrid expressions as program optimization task.
Finally, we demonstrate our approach by generating efficient solvers for the steady-state heat equation with constant and variable coefficients.

%% file: src/grammar.tex
The task of constructing a multigrid solver for a certain problem is typically performed by a human expert with profound knowledge in numerical mathematics.
To automate this task, we first need a way to represent multigrid solvers in a formal language that we can then use to automatically construct different instances on a computer.
The rules of this language must ensure that only valid solver instances can be defined, which means that we can both automatically determine their convergence and runtime behavior.
Additionally, we want to enforce that the generated method works on a hierarchy of grids, which requires the availability of inter-grid operations that allow to obtain approximations of the same operator or grid on a finer or coarser level.
Consider the general system of linear equations defined on a grid with spacing $h$ 
\begin{equation}
	A_{h} u^{h} = f^{h},
	\label{eq:lgs}
\end{equation} 
where $A_{h}$ is the coefficient matrix, $u^{h}$ the unknown and $f^h$ the right-hand side of the system.
Each component of a multigrid solver can be written in the following form:
\begin{equation}
u_{i+1}^h = u_{i}^h + \omega B_{h} (f^h - A_{h} u_{i}^h),\label{eq:stationary-method}
\end{equation}
where $u_{i}^h$ is the approximate solution in iteration $i$, $\omega \in \mathbb{R}$ the relaxation factor and $B_{h}$ an operator defined on the given level with spacing $h$.
For example, with the splitting $A_{h} = D_{h} - U_{h} - L_{h}$, we can define the Jacobi
\begin{equation}
u_{i+1}^h = u_{i}^h + D_{h}^{-1} (f^h - A_{h} u_{i}^h) \label{eq:jacobi}
\end{equation}
and the lexicographical Gauss-Seidel method
\begin{equation}
u_{i+1}^h = u_{i}^h + (D_{h} - L_{h})^{-1} U_{h} (f^h - A_{h} u_{i}^h). \label{eq:lex-gs}
\end{equation}
If we assume the availability of a prolongation operator $P_H$, a restriction operator $R_h$ and an approximation for the inverse of $A_h$ on the coarser grid, a coarse grid correction can be defined as
\begin{equation}
u_{i+1}^h = u_{i}^h + P_H A_H^{-1} R_h(f^h - A_{h} u_{i}^h)\label{eq:cgc}.
\end{equation}
Furthermore, we can substitute $u_{i}^h$ in \eqref{eq:cgc} with \eqref{eq:jacobi} and obtain a two grid with Jacobi pre-smoothing
\begin{equation}
u_{i+1}^h = (u_{i}^h + D_{h}^{-1} (f^h - A_{h} u_{i}^h)) + P_H A_H^{-1} R_h(f^h - A_{h} (u_{i}^h + D_{h}^{-1} (f^h - A_{h} u_{i}^h))). \label{eq:two-grid}
\end{equation}
By repeatedly substituting subexpressions we can automatically construct a single expression for any multigrid solver.
If we take the set of possible substitutions as a basis, we can define a list of rules according to which we can generate such an expression.
We specify these rules in the form of a context-free grammar, which is described in figure \ref{fig:grammar} for classical pointwise smoothers.
\begin{figure}[t]
	\begin{subfigure}[t]{\textwidth}
	\begin{bnf*}
		\bnfprod{$S$} {
			\bnfts{\textsc{iterate}}(\bnfpn{$c^h$}, \bnfsp \bnfts{$\omega$}, \bnfsp \bnfpn{$\mathcal{P}$})
		} \\
		\bnfprod{$c^h$} {
			\bnfts{\textsc{apply}}(\bnfpn{$B_h$}, \bnfsp \bnfpn{$c^h$}) \bnfor \bnfts{\textsc{residual}}(\bnfts{$A_h$}, \bnfsp \bnfpn{$s^h$}) 
		} \\
		\bnfprod{$s^h$} {
			\bnfts{\textsc{iterate}}(\bnfpn{$c^h$}, \bnfsp \bnfts{$\omega$}, \bnfsp \bnfpn{$\mathcal{P}$}) \bnfor 			\bnfts{\textsc{iterate}}(\bnfts{\textsc{coarse-grid-correction}}(\bnfts{$P_{2h}$}, \bnfsp \bnfpn{$c^{2h}$}, \bnfsp \bnfts{$\omega$}), \bnfsp \bnfts{$\omega$}, \bnfsp \bnfes) \bnfor 
			(\bnfts{$u_0^h$}, \bnfsp \bnfts{$f_0^h$},\bnfsp \bnfes, \bnfsp \bnfes)
		} \\
		\bnfprod{$B_h$} {
			\bnfts{\textsc{inverse}}(\bnfts{\textsc{diagonal}}(\bnfts{$A_h$}))
		} \\
		\bnfprod{$c^{2h}$} {
			\bnfts{\textsc{apply}}(\bnfpn{$B_{2h}$}, \bnfsp \bnfpn{$c^{2h}$}) \bnfor
			\bnfts{\textsc{residual}}(\bnfts{$A_{2h}$}, \bnfsp \bnfpn{$s^{2h}$}) \bnfor 
			\bnfts{\textsc{coarse-cycle}}(\bnfts{$A_{2h}$}, \bnfsp \bnfts{$u^{2h}_0$}, \bnfsp \bnfts{\textsc{apply}}(\bnfts{$R_h$}, \bnfsp \bnfpn{$c^h$})) 
		} \\
		\bnfprod{$s^{2h}$} {
			\bnfts{\textsc{iterate}}(\bnfpn{$c^{2h}$}, \bnfsp \bnfts{$\omega$}, \bnfsp \bnfpn{$\mathcal{P}$}) \bnfor
			\bnfts{\textsc{iterate}}(\bnfts{\textsc{apply}}(\bnfts{$P_{4h}$}, \bnfsp \bnfpn{$c^{4h}$}), \bnfsp \bnfts{$\omega$}, \bnfsp \bnfes)
		} \\
		\bnfprod{$B_{2h}$} {
			\bnfts{\textsc{inverse}}(\bnfts{\textsc{diagonal}}(\bnfts{$A_{2h}$}))
		} \\
		\bnfprod{${c}^{4h}$} {
			\bnfts{\textsc{apply}}(\bnfts{$A^{-1}_{4h}$}, \bnfsp \bnfts{\textsc{apply}}(\bnfts{$R_{2h}$}, \bnfsp \bnfpn{$c^{2h}$}))
		} \\
		\bnfprod{$\mathcal{P}$} {
			\bnfts{\textsc{red-black partitioning}} \bnfor \bnfes
		}
	\end{bnf*}
\caption{Syntax}
\label{fig:grammer-syntax}
\end{subfigure}
\begin{subfigure}[b]{\textwidth}
	\begin{algorithmic}
	\Function{iterate}{($u$, $f$, $\delta$, $state$), $\omega$, $\mathcal{P}$}
		\State $\tilde{u} \gets u + \omega \cdot \delta$ with $\mathcal{P}$
		\State return ($\tilde{u}$, $f$, $\lambda$, $state$) 
	\EndFunction
	\Function{apply}{$B$, ($u$, $f$, $\delta$, $state$)}
		\State $\tilde{\delta} \gets B \cdot \delta$  
		\State return ($u$, $f$, $\tilde{\delta}$, $state$)
	\EndFunction
	\Function{residual}{$A$, ($u$, $f$, $\lambda$, $state$)}
		\State $\delta \gets f - A u$
		\State return ($u$, $f$, $\delta$, $state$)
	\EndFunction
	\Function{coarse-cycle}{$A_{H}$, $u^H_0$, ($u^h$, $f^h$, $\delta^H$, $state^h$)}
		\State $u^H \gets u^H_0$ 
		\State $f^H \gets \delta^H$
		\State $\tilde{\delta}^H \gets f^H - A_H u^H_0$ 
		\State $state^H \gets$ ($u^h$, $f^h$, $\delta^H$, $state^h$)
		\State return ($u^H$, $f^H$, $\tilde{\delta}^H$, $state^H$)
	\EndFunction
	\Function{coarse-grid-correction}{$P_H$, ($u^H$, $f^H$, $\delta^H$, $state^H$), $\omega$}
		\State ($u^h$, $f^h$, $\delta^h$, $state^h$) $\gets state^H$
		\State $\tilde{\delta}^h \gets P_H \cdot (u^H + \omega \cdot \delta^H)$
		\State return ($u^h$, $f^h$, $\tilde{\delta}^h$, $state^h$)
	\EndFunction
	\end{algorithmic}
\caption{Semantics}
\label{fig:grammer-semantics}
\end{subfigure}
\caption{A formal grammar for the generation of multigrid solvers.}
\label{fig:grammar}
\end{figure}
Figure \ref{fig:grammer-syntax} contains the production rules while figure \ref{fig:grammer-semantics} describes their semantics.
Each rule defines the set of expressions by which a certain production symbol, denoted by $\langle \cdot \rangle$, can be replaced.
To generate an expression, starting with the symbol $\langle S \rangle$, this process is recursively repeated until the produced expression contains only terminal symbols or the empty string $\lambda$.
The construction of a multigrid solver comprises the recursive generation of cycles on multiple levels.
Consequently, we must be able to create a new system of linear equations on a coarser level, including a new initial solution, right-hand side and coefficient matrix.
%Here we assume that the system matrix $A_h$ is defined on arbitrary levels, alternatively one could employ a Galerkin coarsening.
Moreover, if we decide to finish the computation on a certain level, we must be able to restore the state of the next finer level, i.e. the current solution and right-hand side, when applying the coarse grid correction.
The current state of a multigrid solver on a certain level with grid spacing $h$ is represented as a tuple ($u^h$, $f^h$, $\delta^h$), where $u^h$ represents the current iterate, $f^h$ the right-hand side and $\delta^h$ a correction expression.
To restore the current state on the next finer level, we additionally include a reference $state^h$ to the corresponding state tuple.
According to figure \ref{fig:grammer-syntax}, the construction of a multigrid solver always ends when the tuple ($u^h_0$, $f^h_0$, $\lambda$, $\lambda$) is reached.
This tuple contains the initial solution and right-hand side on the finest level and therefore corresponds to the original system of linear equations that we aim to solve.
Here we have neither yet computed a correction, nor do we need to restore any state and hence both $\delta^h$ and $state^h$ contain the empty string.
In general, our grammar includes three functions that operate on a fixed level.

The function \textsc{iterate} generates a new state tuple from a given one by applying the correction $\delta^h$ to the current iterate $u^h$ using the relaxation factor $\omega$.
If available a partitioning can be included to perform the update in multiple sweeps on subsets of $u^h$ and $\delta^h$, e.g. a red-black Gauss-Seidel iteration.
%Note that the order of evaluation is not considered here, which means that we postpone the decision whether to evaluate a certain subexpression and store its result in a temporary variable to the code generation phase.
The function \textsc{residual} creates a residual expression based on the given state tuple, which is assigned to the newly created symbol $\delta$.
A correction $\delta$ can be transformed with the function \textsc{apply}, which generates a new correction $\tilde{\delta}$ by applying the linear operator $B_h$ to the old one.
For example, the following function applications evaluate to one iteration of damped Jacobi smoothing:
\begin{align*}
& \textsc{iterate}(\textsc{apply}(\textsc{inverse}(\textsc{diagonal}(A_h)),\; \textsc{residual}(A_h,\; (u^h_0,\; f^h,\; \lambda,\; \lambda))), \; 0.7, \; \lambda) \\
& \to \textsc{iterate}(\textsc{apply}(D_h^{-1},\; (u^h_0,\; f^h,\; f^h - A_hu^h_0,\; \lambda)), \; 0.7, \; \lambda) \\
& \to \textsc{iterate}((u^h_0,\; f^h,\; D_h^{-1}(f^h - A_hu^h_0),\; \lambda), \; 0.7, \; \lambda) \\
& \to (u_{0}^h + 0.7 \cdot D_h^{-1} (f^h - A_{h} u_{0}^h), \; f^h, \; \lambda, \; \lambda),
\end{align*}
where $D^{-1}_h$ is the inverse diagonal of $A_h$, as defined above.

Finally, it remains to be shown how one can recursively create a multigrid cycle on the next coarser level and then apply the result of its computation to the current approximate solution.
This is accomplished through the functions \textsc{coarse-cycle} and \textsc{coarse-grid-correction}.
The former expects a state tuple to which the restriction $R_h$ has been already applied.
It then creates a new state on the next coarser level using the initial solution $u_0^H$, the operator $A_H$ and the restricted correction $\delta^H$ as right-hand side $f^H$.
Note that on the coarsest level the resulting system of linear equation can be solved directly, which is denoted by the application of the inverse coarse grid operator.
To enable the reestablishment of the previous state, a reference is stored in $state^H$.
If the computation is finished, the function \textsc{coarse-grid-correction} comes into play.
It first restores the previous state on the next finer level and then computes a coarse grid correction by applying the interpolation operator to the solution computed on the coarser grid, which is then used as a new correction $\tilde{\delta}$ on the finer level.
Again, the following example application demonstrates the semantics of these functions, where we abbreviate \textsc{coarse-grid-correction} with \textsc{cgc}: 
\begin{align*}
	& \textsc{cgc}(P_{2h}, \textsc{coarse-cycle}(A_{2h},\; u_0^{2h}, \; (u^h_0, \; f^h, \; R_h (f^h - A_h u^h_0), \; \lambda)), \; \omega )\\
	& \to \textsc{cgc}(P_{2h}, (u_0^{2h}, \; R_h (f^h - A_h u^h_0), \; R_h (f^h - A_h u^h_0) - A_{2h} u^{2h}_0, \; (u^h_0, \; f^h, \; R_h (f^h - A_h u^h_0), \; \lambda)), \; \omega ) \\
	& \to (u_0^h, \; f^h, \; P_{2h} \cdot (u_0^{2h} + \omega \cdot (R_h (f^h - A_h u^h_0) - A_{2h} u_0^{2h})), \; \lambda)
\end{align*}

With the implementation of these functions we have completed the definition of the syntax and semantics of our formal grammar for the generation of multigrid solvers.
It must be mentioned that this grammar imposes certain restrictions on the structure of the generated solver.
First of all, we only allow the application of operators, either for smoothing, restriction or prolongation, to the current residual.
This is a sufficient constraint for the generation of multigrid solvers for linear system, though nonlinear multigrid methods, such as the full approximation scheme (FAS) \cite{briggs2000multigrid}\cite{trottenberg2000multigrid}, require the restriction and prolongation of the current solution and can therefore not be generated with the grammar presented here.
Furthermore, we assume that the right-hand side $f^h$ is only available on the finest grid.
A full multigrid (FMG) scheme starts on the coarsest grid and hence requires the availability of a right-hand side on each level, which could for instance be introduced as additional terminal symbol.
Besides these obvious restrictions one could arbitrarily loosen the constraints implicitly made within the grammar and enable the combination of coarse grid corrections that originate from different expressions on a finer level.
Even though we can not preclude that it is possible to generate improved multigrid methods without these restrictions, this work only represents a first step towards the automatic generation and optimization of these methods and we do not claim to consider all possible variations, but instead focus on the classical multigrid formulation, as presented in \cite{hackbusch2013multi,briggs2000multigrid,trottenberg2000multigrid}.
Since we have shown how it is possible to generate expressions that uniquely represent different multigrid solvers using the formal grammar defined in figure \ref{fig:grammar}, the remainder of this paper focuses on the evaluation and optimization of the resulting algorithms based on this representation.

%% file: src/optimization.tex
The fundamental requirement for the optimization of an iterative method is to have a way to evaluate both its rate of convergence and performance on the target machine.
As we want to fully automate this process, we must be able to perform all steps from the generation of a solver to its evaluation without requiring any human intervention.
In general, there are two possibilities to automatically evaluate an algorithm.
Assuming there exists a code generator that is able to generate machine instructions from a high-level algorithm, one could first translate it to this representation, then employ the generator to emit an executable and finally run it to evaluate both objectives.
The main disadvantage of this approach is that, depending on the execution time of the code generator, this can be infeasible.
The second possibility is to use predictive models to obtain an approximation for both objectives in significantly less compute time.
%In the following, we will present how one can employ both possibilities to evaluate and optimize the performance of geometric multigrid solvers based on the grammar presented in the last section.
This work focuses on the automatic optimization of geometric multigrid solvers on rectangular grids.
In this case we can represent all matrices as one or multiple stencil codes and there exist models, which we briefly explain in the following, that allow us to predict the quality of a multigrid solver with respect to both objectives.
Although, as it can not be expected that these predictions are always accurate, we still employ code generation to evaluate the best solvers of each optimization run.
\subsection{Convergence estimation}
The quality of an iterative method is first and foremost determined by its rate of convergence, i.e. the speed at which the approximation error approaches machine precision.
One iteration of a multigrid solver can be expressed in the general form of equation \eqref{eq:stationary-method}.
By separating all terms that contain the current iterate $u_i^h$ from the rest of the equation, we obtain the following form:
\begin{equation}\label{eq:multigrid_general_formulation}
u_{i+1}^h = (I_h - \omega B_h A_h) u_i^h + \omega B_h f^h,
\end{equation}
where $I_h$ is the unit matrix.
The \emph{iteration matrix} $M_h$ of the given multigrid solver is then given by
\begin{equation}
M_h = (I_h - \omega B_h A_h).
\end{equation}
The \emph{spectral radius} $\rho$ of this matrix, as defined by
\begin{equation}\label{eq:spectral_radius} 
\rho(M_h) = \max_{1 \leq j \leq n}|\lambda_j(M_h)|,
\end{equation}
where $\lambda_j(M_h)$ are the eigenvalues of $M_h$, is essential for the convergence of the method.
Assume $u_{*}^h$ is the exact solution of the system, the error $e_i^h = u_i^h - u_{*}^h$ in iteration $i$ then satisfies,
\begin{equation}
	e_i^h = M_h^i e_0^h,
\end{equation}
where $M^i_h$ is the $i$th power of $M_h$.
The \emph{convergence factor} of this sequence is the limit
\begin{equation}
\rho = \lim_{i \to \infty} \left( \frac{\norm{e_i^h}}{\norm{e_0^h}} \right)^{1/i},
\end{equation}
which is equal to the spectral radius of the iteration matrix $M_h$ \cite{saad2003iterative}.
In general, the computation of the spectral radius is of complexity $\mathcal{O} (n^3)$ for $M_h \in \mathbb{R}^{n \times n}$.
If we however restrict ourselves to geometric multigrid solvers on rectangular grids, we can employ local Fourier analysis (LFA) to obtain an estimate for $\rho$ \cite{wienands2004practical}.
LFA considers the original problem on an infinite grid while the boundary conditions are neglected.
Recently LFA has been automated through the use of software packages \cite{wienands2004practical,BRFourier2018}.
LFA Lab\footnote{LFA Lab: \url{https://github.com/hrittich/lfa-lab}} is a library for the automatic local Fourier analysis of constant and periodic stencils \cite{BRFourier2018} on rectangular grids.
To automatically estimate the convergence factor of a multigrid solver using this tool, we first need to obtain the iteration matrix.
Using the grammar described in the last section, we always generate expressions of the form of equation \eqref{eq:stationary-method} from which we can extract the iteration matrix using the transformation described above.
Finally, we emit the resulting expression, which represents the iteration matrix of our multigrid solver, in the form of an LFA Lab expression, for which we can automatically estimate the spectral radius.
\subsection{Performance estimation}
A popular yet simple model for estimating the performance of an algorithm on modern computer architectures is the \emph{roofline model} \cite{williams2009roofline}.
Based on the operational intensity of a compute kernel, i.e. the ratio of floating point operations to words loaded from and stored to memory, it gives an estimate for the maximum achievable performance, which is either limited by the memory bandwidth or the compute capabilities of the machine.
The basic roofline formula is given by
\begin{equation}
	P = \min(P_{max}, \; I \cdot b_s),
\end{equation}
where $P$ is the attainable performance, $P_{max}$ the peak performance of the machine, i.e. the maximum achievable amount of floating point operations per second, $I$ the operational intensity of the kernel and $b_s$ the peak memory bandwidth, i.e. the amount of words that can be moved from and to main memory per second.
Within a geometric multigrid solver each kernel either represents a matrix-vector or vector-vector operation, where each vector corresponds to a rectangular grid and each matrix to one or multiple stencils.
If we explicitly represent each operation in the form of a stencil, the computation of the operational intensity is straightforward.
%It must be noted that the roofline model has been extended to obtain more accurate predictions through the use of more sophisticated machine model which we do not consider in this work, such as the Execution-Cache-Memory (ECM) model \cite{hager2016exploring}\cite{treibig2009introducing}.
\subsection{Code Generation and Evaluation}
In order to evaluate a solver that has been evolved within a certain stage of optimization on the target platform, we employ the ExaStencils code generation framework \cite{lengauer2014exastencils}, which was specifically designed for the generation of geometric multigrid implementations that run on parallel and distributed systems. 
To employ the code generation capabilities of this framework, we transform the evolved multigrid expression to an algorithmic representation, which we then emit in the form of a specification in ExaStencils' DSL \cite{schmitt2014exaslang}.
Based on this specification the framework generates a C++ implementation of the solver, including a default application, which we finally run to measure both its total execution time $T$ and defect reduction factor
\begin{equation}
\tilde{\rho}_{i} = \frac{\norm{f^h - A_h u_{i}^h }}{\norm{f^h - A_h u^h_{i-1}}} 
\end{equation}
per iteration $i$ on the target platform. 
We then obtain an approximate for the asymptotic convergence factor
\begin{equation}\label{eq:asymptotic_convergence_factor}
\tilde{\rho} = \left(\prod_{i=1}^{n}\tilde{\rho}_{i} \right)^{1/n},
\end{equation} where $n$ is the number of iterations until convergence \cite{trottenberg2000multigrid}.
\subsection{Optimization}
In case we want to find the optimal geometric multigrid solver for a certain problem, first the question about the size of the search space arises.
With a sufficiently small search space one could attempt to simply enumerate all possible solutions. 
The infeasibility of this approach becomes obvious when looking at the grammar in figure \ref{fig:grammar}.
Assume our goal is to find a multigrid solver that operates on three levels, but the only allowed operation on the coarsest level is the application of a direct solver.
Besides the start symbol $\langle S \rangle$  and the production resulting in the application of a direct solver on the coarsest level, each non-terminal symbol produces at least two alternatives.
Now assume we perform on average twenty productions per level.
This means we must consider more than $2^{40}$ alternatives that must be all evaluated with respect to both objectives, which is already infeasible on a standard desktop computer.
In practice this number will be even larger, especially if we consider more levels.
Furthermore, we need to choose a value for all occurrences of the relaxation parameter $\omega$, which yields an additional continuous optimization problem.
In case the search space is too large to be directly enumerated, a remedy is to use heuristics that aim to search efficiently through the space of possible solutions and still find the global or at least a local optimum.
Evolutionary algorithms are a class of search heuristics inspired by the principle of natural evolution that have been successfully applied to numerous domains \cite{koza2010human}.
All of these methods have in common that they evolve a population of solutions (called individuals) through the iterative application of so-called genetic operators.
Typically each iteration (or generation) of an evolutionary algorithm consists of the following steps:
\begin{enumerate}
	\item \textbf{Selection}: A number of best individuals is selected from the population of the last generation.
	\item \textbf{Crossover}: New individuals are created through the recombination of existing ones.
	\item \textbf{Mutation}: A number of individuals is randomly altered to create new instances.
\end{enumerate}
The order and probability of application of each operation can be varied and different choices have been suggested for different optimization problems \cite{back1997handbook}.
The exact implementation of each genetic operator depends on the class of problem, i.e. the structure of the solution.
Within this work we consider two different optimization problems.
First of all, we want to find the list of productions that, according to the context-free grammar presented in section \ref{sec:grammar}, leads to the optimal multigrid solver.
The class of evolutionary algorithms that evolve expressions according to a context-free grammar are summarized under the term genetic programming \cite{poli08:fieldguide, koza1994genetic}.
To evolve a Pareto front of multigrid solvers with respect to both objectives, we employ strongly typed genetic programming \cite{montana1995strongly} with a non-dominated sorting based selection \cite{coello2007evolutionary}.
Because the computation of the spectral radius significantly slows down for expressions consisting of more than two levels, we split each optimization into multiple runs.
Starting on the three coarsest levels, we perform the optimization assuming that we can obtain the correct solution on the coarsest level.
During this process the value $1$ is chosen for each relaxation factor $\omega$ that occurs within an expression.
After we have evolved a Pareto front of multigrid expressions, for each of those individuals an implementation is generated, which we then evaluate on the target platform by considering an instance of the problem that we aim to solve on the finest level.
For this purpose we adapt the initial solution, boundary condition and right-hand side of the given problem to the grid defined on the current level.
We then choose the individual that leads to the fastest solver with respect to execution time until convergence on the target platform.

In the second step, we turn our attention to the list of relaxation factors in order to improve the convergence of the best individual evolved in the first step, which corresponds to a single-objective continuous optimization problem, that we solve using a covariance matrix adaptation evolution strategy (CMA-ES) \cite{hansen2001completely}.
To evaluate the convergence of a solver with respect to a certain choice of relaxation factors, we reuse the implementation generated at the end of the first optimization step to measure $\tilde{\rho}$, as defined in equation \eqref{eq:asymptotic_convergence_factor}.
Since we only need to adapt a number of real-valued parameters in the respective source file, we just need to recompile the binary while avoiding the high cost of rerunning the ExaStencils code generator.
To summarize, each optimization consists of the following steps:
\begin{enumerate}
	\item Multi-objective optimization using automated local Fourier analyis (LFA) and performance modeling
	\item Evaluation of Pareto-optimal solvers using code generation
	\item Relaxation parameter optimization based on the generated implementation
\end{enumerate}
Finally, the resulting individual is employed as direct solver for the performance estimation and code generation on the next two levels.
We repeat this procedure until the optimization on the finest level is finished.
By recursively deploying the best individual of a run as direct solver for the next run, a single multigrid expression is obtained which operates on the complete range of levels.
We implement our optimization approach in the Python programming language\footnote{EvoStencils: \url{https://github.com/jonas-schmitt/evostencils}} using the framework DEAP \cite{DEAP_JMLR2012} for the implementation of the evolutionary algorithms.

%% file: src/experiments.tex
For our experiments we consider the steady-state heat equation with Dirichlet boundary conditions on a unit square, which is given by
\begin{equation}
    \label{eq:heat-eq}
    \begin{aligned}
        -\divergence \left(a \grad u \right) &= f \quad \text{in} \quad \Omega \,, \\
        u &= g \quad \text{on} \quad \partial \Omega
        \,.
    \end{aligned}
\end{equation}
where $\Omega = (0,1)^d$, $\divergence v: \mathbb{R}^d \to \mathbb{R}$ is the divergence of $v$ and $\grad u: \mathbb{R} \to \mathbb{R}^d$ is the gradient of $u$.
The function $a: \mathbb{R}^d \to \mathbb{R}$ describes the thermal conductivity of the material.
We discretize equation \eqref{eq:heat-eq} using finite differences on a cartesian grid with a step size of $h$ to obtain the system of linear equations $A_h u^h = f^h$.
Our goal is to evolve optimal multigrid methods for solving this system.
For this purpose we consider four different cases which are summarized in figure \ref{figure:test-cases}.
\begin{figure}
	\begin{subfigure}[b]{.49\textwidth}
	\begin{align*}
	f(x,y) & = \pi^2 \cos(\pi x ) - 4 \pi^2 \sin ( 2 \pi y) \\
	a(x,y) & = 1 \\
	g(x,y) & = \cos ( \pi x ) - \sin ( 2 \pi y )
	\end{align*}\caption{2D with constant coefficients}
	\end{subfigure}
	\begin{subfigure}[b]{.49\textwidth}
		\begin{align*}
		f(x,y,z) & = x^2 - 0.5 y^2 - 0.5 z^2 \\
		a(x,y,z) & = 1 \\
		g(x,y,z) & = 0
		\end{align*}\caption{3D with constant coefficients}
	\end{subfigure}
	\begin{subfigure}[b]{.49\textwidth}
	\begin{align*}
	f(x,y) & = 2 \kappa ((x-x^2)+(y-y^2)) \\
	a(x,y) & = e^{\kappa (x-x^2) (y-y^2)} \\
	g(x,y) & = 1 - e^{ (-\kappa) ( x - x^2  ) ( y - y^2  )}
	\end{align*}\caption{2D with variable coefficients}
	\end{subfigure}
	\begin{subfigure}[b]{.49\textwidth}
	\begin{align*}
	& f(x,y,z) = 2\kappa ((x-x^2)(y-y^2) \\ & \quad + (x-x^2)(z-z^2)
	+ (y-y^2)(z-z^2) ) \\
	& a(x,y,z) = e^{k(x-x^2)(y-y^2)(z-z^2)} \\
	& g(x,y,z) = 1 - e ^{ (-\kappa) ( x - x^2  ) (y - y^2) (z - z^2)}
	\end{align*}\caption{3D with variable coefficients}
\end{subfigure}\caption{Overview of the considered test cases.}	\label{figure:test-cases}
\end{figure}
To obtain the same operator on a coarser level, we rediscretize equation \eqref{eq:heat-eq} on a cartesian grid of the appropriate size.
Note that the choice of a constant coefficient function $a(\vec{x}) = 1$ results in Poisson's equation.
To obtain a Pareto front of multigrid expressions, we perform a multi-objective optimization for $100$ generations using genetic programming (GP) with a $(\mu + \lambda)$ ES \cite{beyer2002evolution} with $\mu = \lambda = 1000$, an initial population of $10 \mu$ and the non-dominated sorting procedure presented in \cite{deb2000fast}.
This means that in each generation we create $\lambda$ individuals based on an existing population of size $\mu$ and then select the best $\mu$ individuals for the next generation from the combined set.
The fitness of each individual consists of two objectives, the spectral radius of its iteration matrix, estimated with LFA, and its execution time on the target platform, an Intel Xeon E3-1275 v5 (Skylake) machine with a peak performance of $230.4$ GFLOP/s (four physical cores, 16 DP FLOPs per cycle, $3.6$ GHz clock frequency) and a peak memory bandwidth of $34.1$ GB/s, estimated with the roofline model.
To estimate the spectral radius in the case of variable coefficients, we approximate the coefficient function with a constant stencil at the center of the domain.
Individuals are selected for crossover and mutation using a dominance-based tournament selection as described in \cite{deb2000fast}.
New individuals are created by either crossover with a probability of $0.7$, whereby we employ single-point crossover with a probability of $0.2$ to choose a terminal as crossover point or by mutation, through replacement of a certain subexpression with a new randomly created expression.
To optimize the relaxation factors of a multigrid solver, we employ a CMA-ES \cite{hansen2001completely} with $\lambda = 2 \cdot \lfloor 4 + 3 \ln(n) \rfloor$, where $n$ is the number of relaxation factors, and $200$ generations.
We restrict the set of productions for the generation of operator expressions to
$\langle B_h \rangle \models \textsc{inverse}(\textsc{diagonal}(A_h))$, which means that we only consider pointwise smoothers.
The optimization is performed with a step size of $h = 1/2^{l}$ on each level $l$, whereby we employ a level range of $l \in \left[2,10\right]$ for the 2D and $l \in \left[3,7\right]$ for the 3D cases.
 
\begin{figure}
	\begin{subfigure}{0.49\textwidth}
		\includegraphics[width=\textwidth]{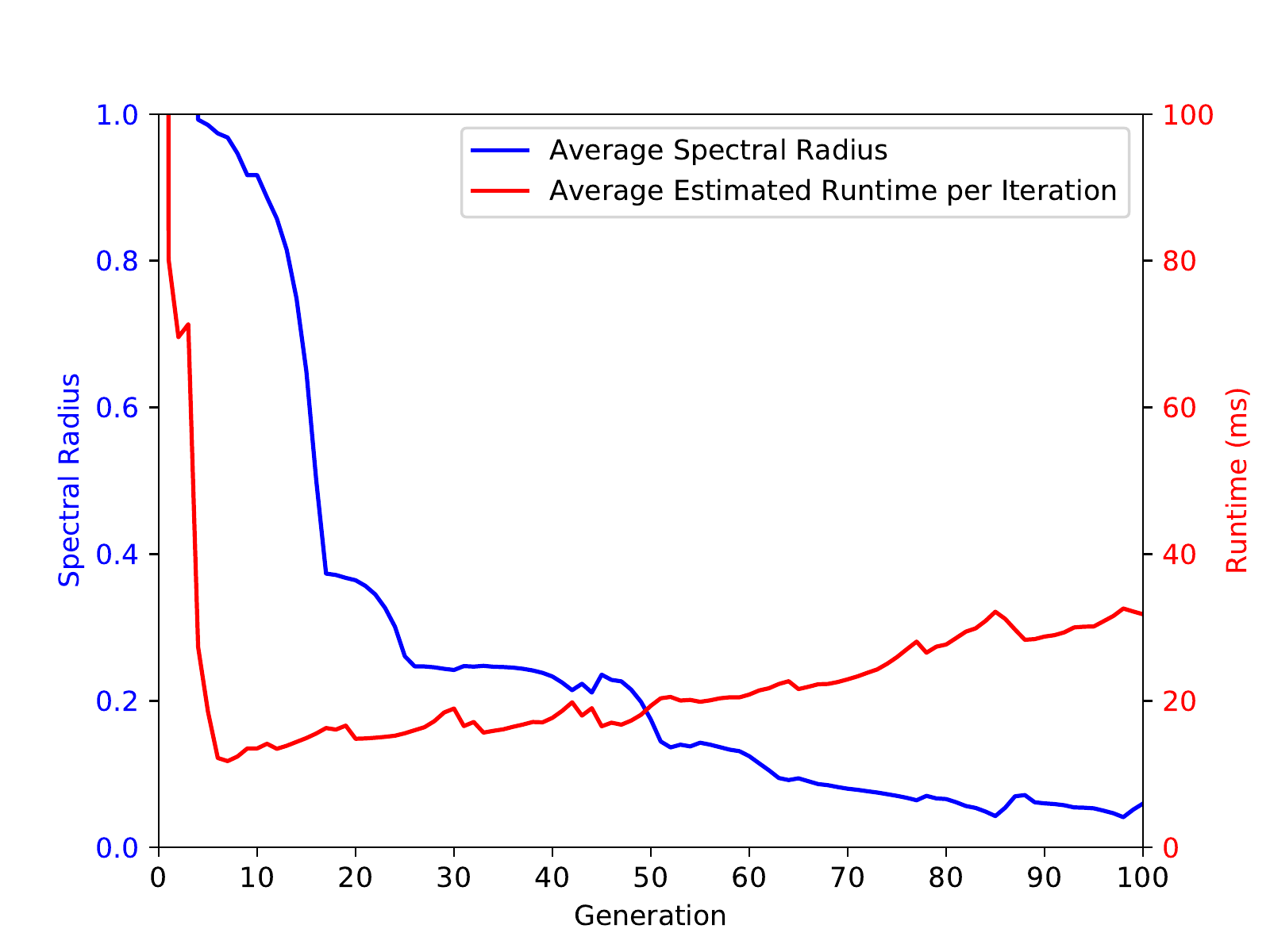}
		\caption{2D with constant coefficients and $l \in \{8, 9, 10\}$}\label{fig:avg-fitness:2D-constant}
	\end{subfigure}
	\begin{subfigure}{0.49\textwidth}
		\includegraphics[width=\textwidth]{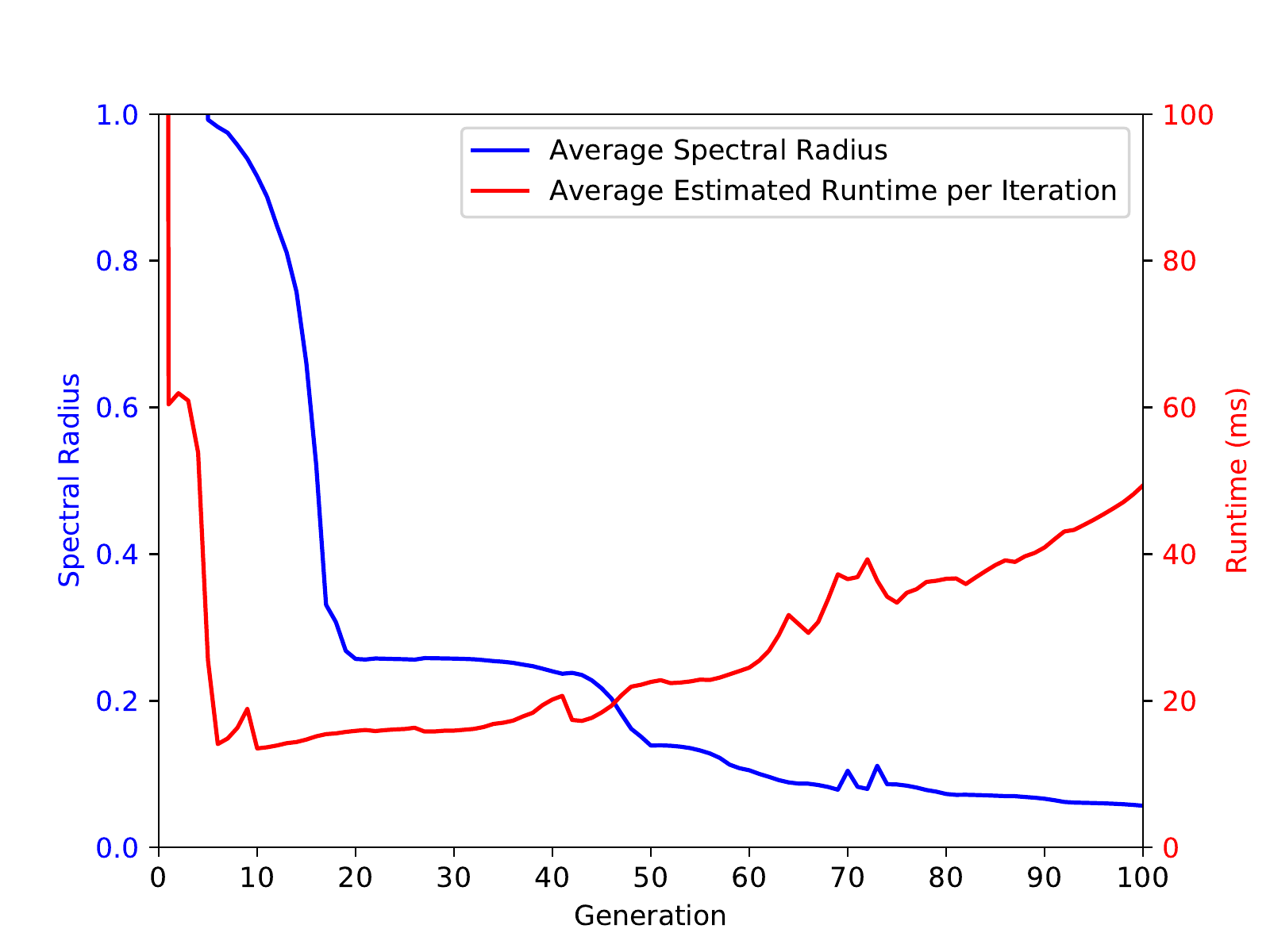}
		\caption{2D with variable coefficients and $l \in \{8, 9, 10\}$}\label{fig:avg-fitness:2D-variable}
	\end{subfigure}
	\begin{subfigure}{0.49\textwidth}
		\includegraphics[width=\textwidth]{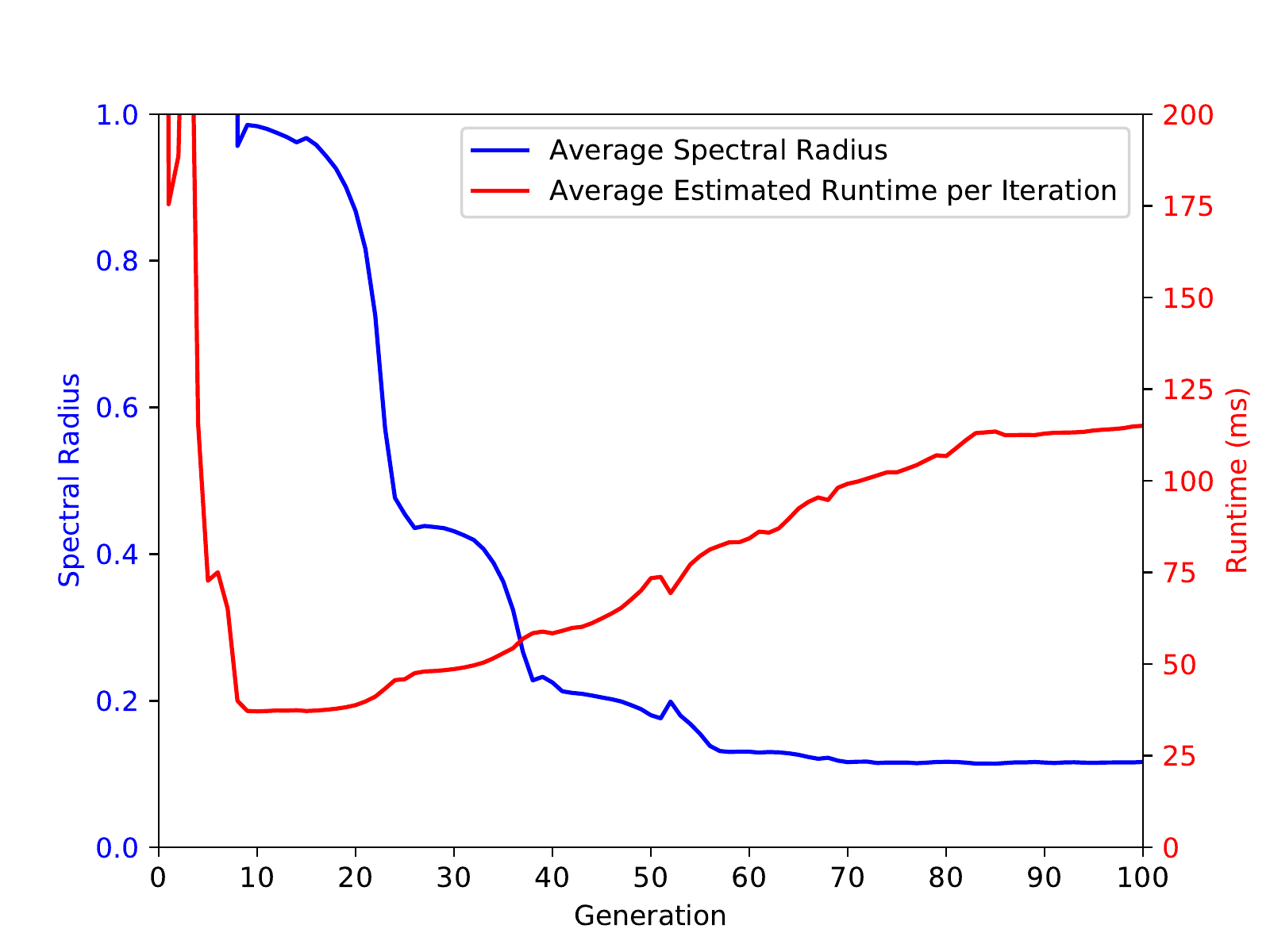}
		\caption{3D with constant coefficients and $l \in \{5, 6, 7\}$}\label{fig:avg-fitness:3D-constant}
	\end{subfigure}
	\begin{subfigure}{0.49\textwidth}
		\includegraphics[width=\textwidth]{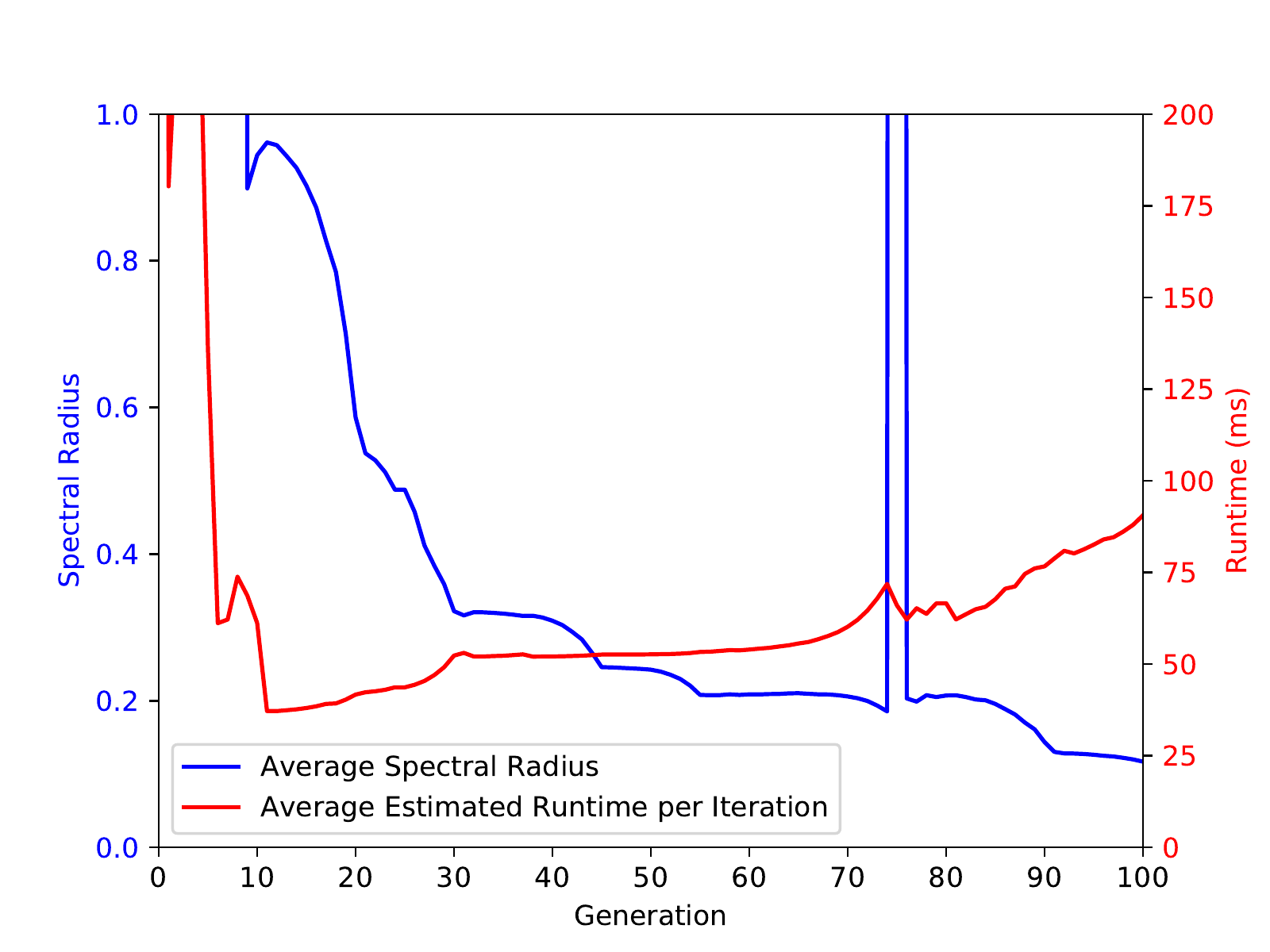}
		\caption{3D with variable coefficients and $l \in \{5, 6, 7\}$}\label{fig:avg-fitness:3D-variable}
	\end{subfigure}
	\caption{Average fitness throughout the optimization.}
	\label{fig:avg-fitness}
\end{figure}

The figures \ref{fig:avg-fitness:2D-constant}, \ref{fig:avg-fitness:2D-variable}, \ref{fig:avg-fitness:3D-constant} and \ref{fig:avg-fitness:3D-variable} each contain a plot of the average value of both objectives in the last GP-optimization of the respective test case, i.e. the optimization on the three finest levels.
In all four cases the evolutionary algorithm is able to decrease the average of both objectives drastically within the first $20 - 30$ generations.
Although in subsequent generations it is not possible to further decrease the average in one objective without entailing an increase in the other one.
This is reflected in the further decrease of the average spectral radius throughout the optimization, while the average execution time per iteration increases slightly.
Since we are first and foremost interested in finding multigrid solvers with fast convergence, we can tolerate a slight increase in the average execution time per iteration if this allows us to find better converging solvers, as the employed selection scheme always preserves those solutions that perform best in both objectives\cite{deb2000fast}.
Since the spectral radius of an iteration matrix does not behave smooth with respect to small perturbations in the structure of a solver, drastic changes in the average can occur throughout an optimization.
Due to the stochastic nature of evolutionary algorithms, for a non-smooth objective function there is always a chance that significantly worse individuals are created through mutation, which is the case in figure \ref{fig:avg-fitness:3D-variable} between the generations 70 and 80.
However, such fluctuations do not affect the general trend of decreasing spectral radii within the population and besides the mentioned deviation, the average of both objectives varies only smoothly throughout each of the four optimization runs, which means that the algorithm possesses a high degree of robustness.
Table \ref{table:convergence_optimization} shows the measured asymptotic convergence factor of the best individual, i.e. the solver with the lowest execution time until convergence, for each test case before and after the convergence optimization using CMA-ES.
\begin{table}\caption{Convergence optimization results on the finest level.}\label{table:convergence_optimization}
	\centering
	\begin{tabular}{c|c|c|c}
		Case & Initial Convergence factor & Optimized convergence factor & Improvement\\
		\hline
		2D with constant coefficients &$0.0740$ & $0.0311$ & $58$ \% \\
		\hline
		2D with variable coefficients &$0.0916$ & $0.0516$ & $44$ \% \\
		\hline
		3D with constant coefficients &$0.0734$ & $0.0190$ & $74$ \% \\
		\hline
		3D with variable coefficients & $0.1409$ & $0.0555$ & $61$ \% \\
		\hline
	\end{tabular}
\end{table}
In all cases a significant improvement can be achieved, which ranges from $44$ to $74$ \%.

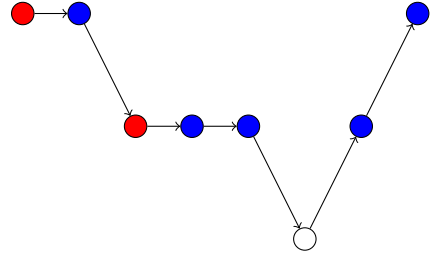
\begin{figure}
\begin{subfigure}[b]{0.395\textwidth}
	\centering
	\begin{algorithmic}[1]
	\State $u^{h} = \textsc{rb-gs}(u^h, f^h, A_h, \omega \approx 1.108)$
	\State $u^{h} =  \textsc{jacobi}(u^h, f^h, A_h, \omega \approx 1.149)$
	\State $\quad f^{2h} = R_h (f^h - A_h u^h)$
	\State $\quad u^{2h} = \textsc{rb-gs}(0, f^{2h}, A_{2h}, \omega \approx 1.031)$
	\State $\quad u^{2h} = \textsc{jacobi}(u^{2h}, f^{2h}, A_{2h}, \omega \approx 0.605)$
	\State $\quad u^{2h} = \textsc{jacobi}(u^{2h}, f^{2h}, A_{2h}, \omega \approx 0.554)$
	\State $\qquad \text{solve} \; A^{4h} u^{4h} = R_{2h} (f^{2h} - A_{2h} u^{2h})$
	\State $\quad u^{2h} = u^{2h} + \omega \cdot P_{4h} u^{4h}$ with $\omega \approx 0.843$
	\State $\quad u^{2h} = \textsc{jacobi}(u^{2h}, f^{2h}, A_{2h}, \omega \approx 0.893)$
	\State $u^{h} = u^{h} + \omega \cdot  P_{2h} u^{2h}$ with $\omega \approx 1.210$
	\State $u^{h} =  \textsc{jacobi}(u^h, f^h, A_h, \omega \approx 1.125)$
\end{algorithmic}
\end{subfigure}
\begin{subfigure}[b]{0.595\textwidth}
	\centering
	\begin{tikzpicture}[nodes=draw]
	\node	(a) at (0,3) [circle,fill=red] {};
	\node	(b) at (0.75,3) [circle,fill=blue] {};
	\node	(c) at (1.5,1.5) [circle,fill=red] {};
	\node	(d) at (2.25,1.5) [circle,fill=blue] {};
	\node	(e) at (3,1.5) [circle,fill=blue] {};
	\node	(f) at (3.75,0) [circle] {};
	\node	(g) at (4.5,1.5) [circle,fill=blue] {};
	\node	(h) at (5.25,3) [circle,fill=blue] {};
	\draw	(a) edge[->] (b)   
	(b) edge[->] (c)
	(c) edge[->] (d)
	(d) edge[->] (e)   
	(e) edge[->] (f)
	(f) edge[->] (g)
	(g) edge[->] (h);
	\end{tikzpicture}
\end{subfigure}
\caption{Evolved solver for the 2D steady-state heat equation with constant coefficients.}\label{fig:algorithm-2D-constant} 
\end{figure}

Finally, to objectively assess the effectiveness of our complete optimization approach, we compare the best solver in each of the four cases with a number of common V- and W-cycles, for each of which we generate a thread-parallel implementation on the target platform using ExaStencils' code generation capabilities.
In all generated implementations we employ statically schedules OpenMP-threads for the parallelization, whereby the number of threads always matches the number of available physical cores.
For compilation on the target machine GCC 8.2.0 is used.
The results of the comparison are listed in the tables \ref{table:solver-comparison:2D-constant}, \ref{table:solver-comparison:2D-variable}, \ref{table:solver-comparison:3D-constant} and \ref{table:solver-comparison:3D-variable}.
\begin{table}\caption{Solver comparison for the steady-state heat  equation.}\label{table:solver-comparison}
		\begin{subtable}{\textwidth}\caption{2D with constant coefficients.}\label{table:solver-comparison:2D-constant}
\centering
\begin{tabular}{c|c|c}
	Method & Convergence factor & Time to convergence \\
	\hline
	V(1,1) & $0.2441$ & $186$ ms \\
	\hline
	V(2,1) & $0.1380$ & $160$ ms \\
	\hline
	V(2,2) & $0.0804$ & $151$ ms \\
	\hline
	V(3,3) & $0.0451$ & $155$ ms \\
	\hline
	W(2,2) & $0.0589$ & $2320$ ms\\
	\hline
	W(3,3) & $0.0320$ & $1867$ ms \\
	\hline
	EvoStencils & $\mathbf{0.0311}$ & $\mathbf{111}$ \textbf{ms} \\
\end{tabular}
		\end{subtable}
		\begin{subtable}{\textwidth}\caption{2D with variable coefficients}\label{table:solver-comparison:2D-variable}
\centering
\begin{tabular}{c|c|c}
	Method & Convergence factor & Time to convergence \\
	\hline
	V(1,1) & $0.3028$ & $2741$ ms \\
	\hline
	V(2,1) & $0.2125$ & $2672$ ms \\
	\hline
	V(2,2) & $0.1426$ & $2538$ ms \\
	\hline
	V(3,3) & $0.1104$ & $3246$ ms \\
	\hline
	W(2,2) & $0.0486$ & $4606$ ms \\
	\hline
	W(3,3) & $\mathbf{0.0305}$ & $4658$ ms \\
	\hline
	EvoStencils & $0.0516$ & $\mathbf{1868}$ \textbf{ms} \\
\end{tabular}
		\end{subtable}
				\begin{subtable}{\textwidth}\caption{3D with constant coefficients}\label{table:solver-comparison:3D-constant}
\centering
\begin{tabular}{c|c|c}
	Method & Convergence factor & Time to convergence \\
	\hline
	V(1,1) & $0.3962$ & $548$ ms \\
	\hline
	V(2,1) & $0.2574$ & $469$ ms \\
	\hline
	V(2,2) & $0.1679$ & $409$ ms \\
	\hline
	V(3,3) & $0.0927$ & $411$ ms \\
	\hline
	W(2,2) & $0.1435$ & $653$ ms \\
	\hline
	W(3,3) & $0.0670$ & $582$ ms \\
	\hline
	EvoStencils & $\mathbf{0.0190}$ & $\mathbf{226}$ \textbf{ms} \\
\end{tabular}
				\end{subtable}
			\begin{subtable}{\textwidth}\caption{3D with variable coefficients}\label{table:solver-comparison:3D-variable}
\centering
\begin{tabular}{c|c|c}
	Method & Convergence factor & Time to convergence \\
	\hline
	V(1,1) & $0.3966$ & $9095$ ms \\
	\hline
	V(2,1) & $0.2676$ & $8382$ ms \\
	\hline
	V(2,2) & $0.1729$ & $7792$ ms \\
	\hline
	V(3,3) & $0.0973$ & $7500$ ms \\
	\hline
	W(2,2) & $0.1145$ & $7381$ ms \\
	\hline
	W(3,3) & $\mathbf{0.0502}$ & $7279$ ms \\
	\hline
	EvoStencils & $0.0555$ & $\mathbf{4806}$ \textbf{ms} \\
\end{tabular}
			\end{subtable}
			\end{table}
In all four cases the evolved solvers outperform the tested V- and W-cycles by a significant margin, with reductions in execution time until convergence that range from $26$ to $45$ \% compared to the fastest method in each case.
Moreover, all evolved cycles achieve convergence factors comparable to those of W-cycles, even though all of them can be classified as V-cycle.
For a better interpretability of these results, we depict the algorithmic structure of the evolved multigrid methods in the figures \ref{fig:algorithm-2D-constant}, \ref{fig:algorithm-2D-variable}, \ref{fig:algorithm-3D-constant} and \ref{fig:algorithm-3D-variable}
\begin{figure}
	\begin{subfigure}[b]{0.395\textwidth}
		\centering
		\begin{algorithmic}[1]
			\State $u^{h} = \textsc{rb-gs}(u^h, f^h, A_h, \omega \approx 1.012)$
			\State $u^{h} = \textsc{jacobi}(u^h, f^h, A_h, \omega \approx 1.164)$
			\State $\quad f^{2h} = R_h (f^h - A_h u^h)$
			\State $\quad u^{2h} = \textsc{jacobi}(u^{2h}, f^{2h}, A_{2h}, \omega \approx 0.976)$
			\State $\quad u^{2h} = \textsc{jacobi}(u^{2h}, f^{2h}, A_{2h}, \omega \approx 1.005)$
			\State $\quad u^{2h} = \textsc{jacobi}(u^{2h}, f^{2h}, A_{2h}, \omega \approx 0.865)$
			\State $\qquad \text{solve} \; A^{4h} u^{4h} = R_{2h} (f^{2h} - A_{2h} u^{2h})$
			\State $\quad u^{2h} = u^{2h} + \omega \cdot P_{4h} u^{4h}$ with $\omega \approx 0.855$
			\State $\quad u^{2h} = \textsc{jacobi}(u^{2h}, f^{2h}, A_{2h}, \omega \approx 1.002)$
			\State $u^{h} = u^{h} + \omega \cdot  P_{2h} u^{2h}$ with $\omega \approx 1.169$
			\State $u^{h} =  \textsc{jacobi}(u^h, f^h, A_h, \omega \approx 1.078)$
		\end{algorithmic}
	\end{subfigure}
	\begin{subfigure}[b]{0.595\textwidth}
		\centering
		\begin{tikzpicture}[nodes=draw]
		\node	(a) at (0,3) [circle,fill=red] {};
		\node	(b) at (0.75,3) [circle,fill=blue] {};
		\node	(c) at (1.5,1.5) [circle,fill=blue] {};
		\node	(d) at (2.25,1.5) [circle,fill=blue] {};
		\node	(e) at (3,1.5) [circle,fill=blue] {};
		\node	(f) at (3.75,0) [circle] {};
		\node	(g) at (4.5,1.5) [circle,fill=blue] {};
		\node	(h) at (5.25,3) [circle,fill=blue] {};
		\draw	
		(a) edge[->] (b)   
		(b) edge[->] (c)
		(c) edge[->] (d)
		(d) edge[->] (e)   
		(e) edge[->] (f)
		(f) edge[->] (g)
		(g) edge[->] (h)
		;
		\end{tikzpicture}
	\end{subfigure}
	\caption{Evolved solver for the 2D steady-state heat equation with variable coefficients.}\label{fig:algorithm-2D-variable} 
\end{figure}
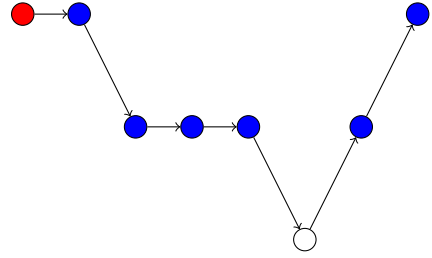
\begin{figure}
	\begin{subfigure}[b]{0.395\textwidth}
		\centering
		\begin{algorithmic}[1]
			\State $u^{h} = \textsc{rb-gs}(u^h, f^h, A_h, \omega \approx 1.545)$
			\State $u^{h} = \textsc{jacobi}(u^h, f^h, A_h, \omega \approx 1.174)$
			\State $u^{h} = \textsc{jacobi}(u^h, f^h, A_h, \omega \approx 1.127)$
			\State $u^{h} = \textsc{jacobi}(u^h, f^h, A_h, \omega \approx 0.745)$
			\State $\quad f^{2h} = R_h (f^h - A_h u^h)$
			\State $\quad u^{2h} = \textsc{rb-gs}(0, f^{2h}, A_{2h}, \omega \approx 0.319)$
			\State $\quad u^{2h} = \textsc{jacobi}(u^{2h}, f^{2h}, A_{2h}, \omega \approx 0.189)$
			\State $\qquad \text{solve} \; A^{4h} u^{4h} = R_{2h} (f^{2h} - A_{2h} u^{2h})$
			\State $\quad u^{2h} = u^{2h} + \omega \cdot P_{4h} u^{4h}$ with $\omega \approx 0.797$
			\State $\quad u^{2h} = \textsc{jacobi}(u^{2h}, f^{2h}, A_{2h}, \omega \approx 1.013)$
			\State $\quad u^{2h} = \textsc{jacobi}(u^{2h}, f^{2h}, A_{2h}, \omega \approx 0.756)$
			\State $\quad u^{2h} = \textsc{jacobi}(u^{2h}, f^{2h}, A_{2h}, \omega \approx 0.689)$
			\State $\quad u^{2h} = \textsc{jacobi}(u^{2h}, f^{2h}, A_{2h}, \omega \approx 1.630)$
			\State $u^{h} = u^{h} + \omega \cdot  P_{2h} u^{2h}$ with $\omega \approx 1.271$
			\State $u^{h} =  \textsc{jacobi}(u^h, f^h, A_h, \omega \approx 1.202)$
		\end{algorithmic}
	\end{subfigure}
	\begin{subfigure}[b]{0.595\textwidth}
		\centering
		\begin{tikzpicture}[nodes=draw]
		\node	(a) at (0,3) [circle,fill=red] {};
		\node	(b) at (0.6,3) [circle,fill=blue] {};
		\node	(c) at (1.2,3) [circle,fill=blue] {};
		\node	(d) at (1.8,3) [circle,fill=blue] {};
		\node	(e) at (2.4,1.5) [circle,fill=red] {};
		\node	(f) at (3.0,1.5) [circle,fill=blue] {};
		\node	(g) at (3.6,0) [circle] {};
		\node	(h) at (4.2,1.5) [circle,fill=blue] {};
		\node	(i) at (4.8,1.5) [circle,fill=blue] {};
		\node	(j) at (5.4,1.5) [circle,fill=blue] {};
		\node	(k) at (6,1.5) [circle,fill=blue] {};
		\node	(l) at (6.6,3) [circle,fill=blue] {};
		\draw	
		(a) edge[->] (b)   
		(b) edge[->] (c)
		(c) edge[->] (d)
		(d) edge[->] (e)   
		(e) edge[->] (f)
		(f) edge[->] (g)
		(g) edge[->] (h)
		(h) edge[->] (i)
		(i) edge[->] (j)
		(j) edge[->] (k)
		(k) edge[->] (l);
		\end{tikzpicture}
	\end{subfigure}
	\caption{Evolved solver for the 3D steady-state heat equation with constant coefficients.}\label{fig:algorithm-3D-constant} 
\end{figure}
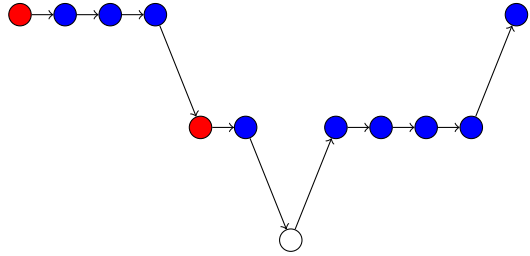
\begin{figure}
	\begin{subfigure}[b]{0.395\textwidth}
		\centering
		\begin{algorithmic}[1]
			\State $u^{h} = \textsc{rb-gs}(u^h, f^h, A_h, \omega \approx 1.093)$
			\State $u^{h} =  \textsc{jacobi}(u^h, f^h, A_h, \omega \approx 1.174)$
			\State $\quad f^{2h} = R_h (f^h - A_h u^h)$
			\State $\quad u^{2h} = \textsc{jacobi}(0, f^{2h}, A_{2h}, \omega \approx 0.534)$
			\State $\quad u^{2h} = \textsc{jacobi}(u^{2h}, f^{2h}, A_{2h}, \omega \approx 0.513)$
			\State $\quad u^{2h} = \textsc{jacobi}(u^{2h}, f^{2h}, A_{2h}, \omega \approx 0.489)$
			\State $\qquad \text{solve} \; A^{4h} u^{4h} = R_{2h} (f^{2h} - A_{2h} u^{2h})$
			\State $\quad u^{2h} = u^{2h} + \omega \cdot P_{4h} u^{4h}$ with $\omega \approx 0.735$
			\State $\quad u^{2h} = \textsc{jacobi}(u^{2h}, f^{2h}, A_{2h}, \omega \approx 1.499)$
			\State $u^{h} = u^{h} + \omega \cdot  P_{2h} u^{2h}$ with $\omega \approx 1.365$
			\State $u^{h} =  \textsc{jacobi}(u^h, f^h, A_h, \omega \approx 1.291)$
		\end{algorithmic}
	\end{subfigure}
	\begin{subfigure}[b]{0.595\textwidth}
		\centering
		\begin{tikzpicture}[nodes=draw]
		\node	(a) at (0,3) [circle,fill=red] {};
		\node	(b) at (0.75,3) [circle,fill=blue] {};
		\node	(c) at (1.5,1.5) [circle,fill=blue] {};
		\node	(d) at (2.25,1.5) [circle,fill=blue] {};
		\node	(e) at (3,1.5) [circle,fill=blue] {};
		\node	(f) at (3.75,0) [circle] {};
		\node	(g) at (4.5,1.5) [circle,fill=blue] {};
		\node	(h) at (5.25,3) [circle,fill=blue] {};
		\draw	(a) edge[->] (b)   
		(b) edge[->] (c)
		(c) edge[->] (d)
		(d) edge[->] (e)   
		(e) edge[->] (f)
		(f) edge[->] (g)
		(g) edge[->] (h);
		\end{tikzpicture}
	\end{subfigure}
	\caption{Evolved solver for the 3D steady-state heat equation with variable coefficients.}\label{fig:algorithm-3D-variable} 
\end{figure}
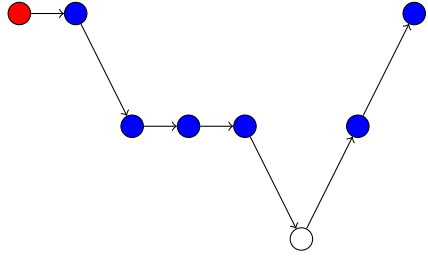
Each of those figures contains both an algorithmic as well as a graphic representation of the evolved multigrid solver on the three finest levels.
In the latter each red node corresponds to one step of red-black Gauss-Seidel smoothing, each blue node to one step of Jacobi smoothing and each white node to the exact solution of the system on the given level.
While the order of operations is displayed from left to right, restriction corresponds to a descend from top to bottom and prolongation to an ascend in the opposite direction.
Even though the evolved methods are structurally different, they exhibit a number of common characteristics.
For the 2D cases the evolved methods are almost identical for constant and variable coefficients in terms of their algorithmic structure, besides the use of a different method as first pre-smoothing step on the second finest level. 
Apart from that, all methods can be characterized as V-Cycles, which means that they only solve once exactly on the coarsest grid.
Furthermore, all of them employ a combination of red-black Gauss-Seidel (RBGS) and Jacobi smoothing with different relaxation factors and a different number of smoothing steps on each level, whereby the latter is typically higher for the coarser levels.
For pre-smoothing often one step of RBGS is used, either solely or followed by a varying number of Jacobi steps.
In contrast, for post-smoothing we observe an exclusive use of under- and overrelaxed Jacobi iterations.
The occurrence of these patterns is especially remarkable in front of the background that we have not prescribed the application of RBGS or Jacobi iterations, but only restricted the function \textsc{apply} within our grammar (see figure \ref{fig:grammar}) to the application of the inverse diagonal of the operator $A_h$ to the current residual, whereby we allow a red-black partitioned computation when obtaining the new iterate.
Even though this obviously includes Jacobi- and partitioned Gauss-Seidel-like methods, it also allows the application of the inverse diagonal zero or multiple times before computing a new iterate.
In the former case the residual is directly added to the current iterate, which corresponds to a modified Richardson iteration.
Guided by the minimization of both objectives, the spectral radius of the resulting iteration matrix and the estimated combined execution time of all applied operations, the evolutionary algorithm nevertheless learns that it is a good idea to apply the inverse diagonal of $A_h$ only once to the current iterate and then repeat this procedure, with different partitioning and relaxation factors.

%% file: src/conclusion.tex
In this work we have presented a novel approach for the automatic optimization of geometric multigrid methods based on a tailored context-free grammar for the generation of multigrid solvers and the use of evolutionary algorithms guided by a model-based prediction for the convergence and compute performance.
Even though we have demonstrated that our approach in principle works and is able to evolve competitive solver for the steady-state heat equation starting from a randomly initialized population of solutions, there is still room for improvement and extensions.
Instead of considering a case that is well researched and for which efficient solution methods are available, a more challenging task would be the solution of partial differential equations where a robust and efficient geometric multigrid solver has not been developed, which is for instance the case for many nonlinear PDEs.
We furthermore want to improve the performance of our evolutionary algorithm by incorporating domain knowledge, in the form of individuals that are known to be efficient in solving a certain problem from long-standing research and analysis in the field of numerical mathematics and which can be directly included into the population to guide the evolution towards promising areas of the search space.
We also aim to improve the accuracy of our model-based compute performance prediction through the use of a more sophisticated machine model such as the Execution-Cache-Memory (ECM) model \cite{hager2016exploring}.
Finally, one could consider different algorithms for the optimization of programs generated by our multigrid grammar, such as reinforcement learning \cite{sutton2018reinforcement} or Monte Carlo tree search \cite{kocsis2006bandit} based techniques.